\documentclass[10pt,reqno]{amsart}
\usepackage{amsmath}
\usepackage{amssymb}
\usepackage{amsthm}

\newlength{\defbaselineskip}
\newcommand{\setlinespacing}[1]%
           {\setlength{\baselineskip}{#1 \defbaselineskip}}

\numberwithin{equation}{section}

\newtheorem{thm}{Theorem}[section]

\theoremstyle{definition}

\theoremstyle{remark}

\numberwithin{equation}{section}

\begin{document}

\title[Unique continuation for fractional Schr\"odinger operators]
{Unique continuation for fractional Schr\"odinger operators in three and higher dimensions}

\author{Ihyeok Seo}


\address{Department of Mathematics, Sungkyunkwan University, Suwon 440-746, Republic of Korea}
\email{ihseo@skku.edu}

\maketitle

\begin{abstract}
We prove the unique continuation property for the differential inequality $|(-\Delta)^{\alpha/2}u|\leq|V(x)u|$,
where $0<\alpha<n$ and $V\in L_{\textrm{loc}}^{n/\alpha,\infty}(\mathbb{R}^n)$, $n\geq3$.
\end{abstract}

\section{Introduction}

In this note we are concerned with the unique continuation property for solutions of
the differential inequality
\begin{equation}\label{equ}
|(-\Delta)^{\alpha/2}u|\leq|V(x)u|,\quad x\in\mathbb{R}^n,\quad n\geq2,
\end{equation}
where $(-\Delta)^{\alpha/2}$, $0<\alpha<n$, is defined
by means of the Fourier transform $\mathcal{F}f$ $(=\widehat{f}\,)$:
$$\mathcal{F}[(-\Delta)^{\alpha/2}f](\xi)=|\xi|^\alpha\widehat{f}(\xi).$$
In particular, the equation $((-\Delta)^{\alpha/2}+V(x))u=0$
has attracted interest from quantum mechanics in the case $1<\alpha<2$
as well as the case $\alpha=2$.
Recently, by generalizing the Feynman path integral to the L\'{e}vy one,
Laskin~\cite{La} introduced the fractional quantum mechanics in which it is conjectured that
physical realizations may be limited to $1<\alpha<2$, where averaged quantities are finite, and
the fractional Schr\"odinger operator $(-\Delta)^{\alpha/2}+V(x)$ plays a central role.
Of course, the case $\alpha=2$ becomes equivalent to an ordinary quantum mechanics.

The unique continuation property means that a solution of ~\eqref{equ}
which vanishes in an open subset of $\mathbb{R}^n$ must vanish identically.
In the case of $\alpha=2$, Jerison and Kenig~\cite{JK} proved the property
for $V\in L_{\textrm{loc}}^{n/2}$, $n\geq3$.
An extension to $L_{\textrm{loc}}^{n/2,\infty}$ was obtained by Stein~\cite{St2}
with small norm in the sense that
$$\sup_{a\in\mathbb{R}^n}\lim_{r\rightarrow0}\|V\|_{L^{n/2,\infty}(B(a,r))}$$
is sufficiently small.
Here, $B(a,r)$ denotes the ball of radius $r>0$ centered at $a\in\mathbb{R}^n$.
These results later turn out to be optimal in the context of $L^p$ spaces (\cite{KN,KT}).

On the other hand, the results when $\alpha\neq2$ are rather scarce.
Laba ~\cite{L} considered the higher orders where $\alpha/2$ are integers,
and obtained the property for $V\in L_{\textrm{loc}}^{n/\alpha}$.
Recently, there was an attempt ~\cite{S2} to handle the non-integer orders when $n-1\leq\alpha<n$, $n\geq2$,
from which it turns out that the condition $V\in L^p$, $p>n/\alpha$, is sufficient to have the property.
Hence this particularly gives a unique continuation result for the fractional Schr\"odinger operator
in the full range $1<\alpha<2$ when $n=2$.
Our aim here is to fill the gap, $0<\alpha<n-1$, for $n\geq3$,
which allows us to have the unique continuation for the fractional Schr\"odinger operator
when $n\geq3$ with the full range of $\alpha$.

\begin{thm}\label{thm}
Let $n\geq3$ and $0<\alpha<n$.
Assume that $V\in L_{\textrm{loc}}^{n/\alpha,\infty}$ and $u$ is a nonzero solution of ~\eqref{equ} such that
\begin{equation}\label{cond}
u\in L^1\cap L^{p,q}\quad\text{and}\quad (-\Delta)^{\alpha/2}u\in L^q,
\end{equation}
where $p=2n/(n-\alpha)$ and $q=2n/(n+\alpha)$.
Then it cannot vanish in any non-empty open subset of $\mathbb{R}^n$ if
\begin{equation}\label{small}
\sup_{a\in\mathbb{R}^n}\lim_{r\rightarrow0}\|V\|_{L^{n/\alpha,\infty}(B(a,r))}
\end{equation}
is sufficiently small.
Here, $L^{p,q}$ denotes the usual Lorentz space.
\end{thm}

\noindent\textit{Remarks.} \textit{(a)}
The smallness condition ~\eqref{small} is trivially satisfied
for $V\in L_{\textrm{loc}}^{n/\alpha}$ because $L_{\textrm{loc}}^{n/\alpha}\subset L_{\textrm{loc}}^{n/\alpha,\infty}$.
Hence the above theorem can be seen as natural extensions to ~\eqref{equ} of the results obtained in~\cite{JK,St2}
for the Schr\"odinger operator ($\alpha=2$).
As an immediate consequence of the theorem, the same result also holds for the stationary equation
$$((-\Delta)^{\alpha/2}+V(x))u=Eu,\quad E\in\mathbb{C},$$
because $(-\Delta)^{\alpha/2}u=(E-V(x))u$ and the condition~\eqref{small} is trivially satisfied for the constant $E$.

\smallskip

\noindent\textit{(b)}
The index $n/\alpha$ is quite natural, in view of the standard rescaling:
$u_\varepsilon(x)=u(\varepsilon x)$ takes the equation $(-\Delta)^{\alpha/2}u=Vu$
into $(-\Delta)^{\alpha/2}u_\varepsilon=V_\varepsilon u_\varepsilon,$
where $V_\varepsilon(x)=\varepsilon^\alpha V(\varepsilon x)$.
So, $\|V_\varepsilon\|_{L^{p,\infty}}=\varepsilon^{\alpha-n/p}\|V\|_{L^{p,\infty}}$.
Hence the $L^{p,\infty}$ norm of $V_\varepsilon$ is independent of $\varepsilon$
precisely when $p=n/\alpha$.

\smallskip

\noindent\textit{(c)}
When $\alpha=n$ in~\eqref{equ}, there are some unique continuation results with
$V\in L_{\textrm{loc}}^p$, $p>1$. (See ~\cite{JK} and ~\cite{S} for $\alpha=2$ and $\alpha=2m$ ($m\in\mathbb{N}$), respectively.)

\

\noindent\textbf{Acknowledgment.}
I am very grateful to Luis Escauriaza for bringing my attention to the papers~\cite{JK,St2}
and for helpful suggestions and discussions.

\section{Proof of the theorem}
From now on, we will use the letter $C$ to denote a constant that may be different at each occurrence.

Without loss of generality, we need to prove that $u$ must vanish identically if
it vanishes in a sufficiently small neighborhood of zero.

Our proof is based on the following Carleman estimate which will be shown below:
If $f\in C_0^\infty(\mathbb{R}^n\setminus\{0\})$ and $(-\Delta)^{\alpha/2}f\in C_0^\infty(\mathbb{R}^n\setminus\{0\})$,
then there is a constant $C$ depending only on $\delta_t:=\min_{k\in\mathbb{Z}}|t-k|$ and $n$
such that for $t\not\in\mathbb{Z}$ with $\delta_t<n-\alpha$
\begin{equation}\label{ineq}
\big\||x|^{-t-n/p}f\big\|_{L^{p,q}}\leq C\big\||x|^{-t+\alpha-n/q}(-\Delta)^{\alpha/2}f\big\|_{L^q},
\end{equation}
where $p,q$ are given as in the theorem (i.e., $1/p+1/q=1$ and $1/q-1/p=\alpha/n$).

Indeed, since we are assuming that $u\in L^1\cap L^{p,q}$ and $(-\Delta)^{\alpha/2}u\in L^q$ vanish near zero
(see~\eqref{cond},~\eqref{equ}),
from ~\eqref{ineq} (with a standard limiting argument involving a $C_0^\infty$ approximate identity), we see that
\begin{equation*}
\big\||x|^{-t-n/p}u\big\|_{L^{p,q}}\leq C\big\||x|^{-t+\alpha-n/q}(-\Delta)^{\alpha/2}u\big\|_{L^q}.
\end{equation*}
Hence,
\begin{align*}
\big\||x|^{-t-n/p}u\big\|_{L^{p,q}(B(0,r))}&\leq C\big\||x|^{-t+\alpha-n/q}Vu\big\|_{L^q(B(0,r))}\\
&+C\big\||x|^{-t+\alpha-n/q}(-\Delta)^{\alpha/2}u\big\|_{L^q(\mathbb{R}^n\setminus B(0,r))}.
\end{align*}
The first term on the right-hand side can be absorbed into the left-hand side as follows:
\begin{align*}
C\big\||x|^{-t+\alpha-n/q}Vu\big\|_{L^q(B(0,r))}&\leq
C\|V\|_{L^{n/\alpha,\infty}(B(0,r))}\big\||x|^{-t+\alpha-n/q}u\big\|_{L^{p,q}(B(0,r))}\\
&\leq \frac12\big\||x|^{-t-n/p}u\big\|_{L^{p,q}(B(0,r))}
\end{align*}
if we choose $r$ small enough (see~\eqref{small}).
Here, recall that $\alpha-n/q=-n/p$, and $\||x|^{-t-n/p}u\|_{L^{p,q}(B(0,r))}$ is finite since $u\in L^{p,q}$ vanishes near zero.
So, we get
$$\|(r/|x|)^{t+n/p}u\|_{L^{p,q}(B(0,r))}
\leq2C\|(-\Delta)^{\alpha/2}u\|_{L^q(\mathbb{R}^n\setminus B(0,r))}<\infty.$$
Now, we choose a sequence $\{t_i\}$ of values of $t$ tending to infinity
such that $\delta_{t_i}$ is independent of $i\in\mathbb{N}$.
Then, by letting $i\rightarrow\infty$, we see that $u=0$ on $B(0,r)$,
which implies $u\equiv0$ by a standard connectedness argument.

\begin{proof}[Proof of ~\eqref{ineq}]
We will show~\eqref{ineq} using Stein's complex interpolation, as in~\cite{St2},
on an analytic family of operators $T_z$ defined by
$$T_zg(x)=\int_{\mathbb{R}^n}K_z(x,y)g(y)|y|^{-n}dy,$$
where $K_z(x,y)=H_z(x,y)/\Gamma(n/2-z/2)$ with
$$H_z(x,y)=|x|^{-t}|y|^{n+t-z}c_z\bigg(|x-y|^{-n+z}-\sum_{j=0}^{m-1}\frac1{j!}\Big(\frac{\partial}{\partial s}\Big)^j|sx-y|^{-n+z}|_{s=0}\bigg).$$
Note that for $f\in C_0^\infty(\mathbb{R}^n\setminus\{0\})$
\begin{equation}\label{22}
T_\alpha(|y|^{-t+\alpha}(-\Delta)^{\alpha/2}f(y))(x)=|x|^{-t}f(x)/\Gamma(n/2-\alpha/2)
\end{equation}
(see Lemma 2.1 in~\cite{S2}).

Let $m$ be a fixed positive integer such that $m-1<t<m$, and recall the following two estimates for the cases of $\textrm{Re}\,z=0$ (Lemma 2.3 in \cite{JK}) and $n-1<\textrm{Re}\,z<n-\delta_t$ (Lemma 4 in \cite{St2}):
There is a constant $C$ depending only on $\delta_t$ and $n$ such that
\begin{equation}\label{JK}
\|T_{i\gamma}g\|_{L^2(dx/|x|^n)}\leq Ce^{c|\gamma|}\|g\|_{L^2(dx/|x|^n)},\quad \gamma\in\mathbb{R},
\end{equation}
and
\begin{equation}\label{St}
\|T_z g\|_{L^r(dx/|x|^n)}\leq Ce^{c|\gamma|}\|g\|_{L^s(dx/|x|^n)},\quad \gamma=\textrm{Im}\,z\in\mathbb{R},
\end{equation}
where $n-1<\beta=\textrm{Re}\,z<n-\delta_t$, $1/s-1/r=\beta/n$ and $1<s<n/\beta$.

We first consider the case where $n-1<\alpha<n$.
Note that we can choose $\beta$ so that $\alpha<\beta<n-\delta_t$, since we are assuming $\delta_t<n-\alpha$.
Hence, by Stein's complex interpolation (\cite{St}) between ~\eqref{JK} and ~\eqref{St}, we see that
\begin{equation}\label{11}
\|T_\alpha g\|_{L^r(dx/|x|^n)}\leq C\|g\|_{L^s(dx/|x|^n)},
\end{equation}
where $1/s-1/r=\alpha/n$ and $1<s<n/\alpha$.
From this and ~\eqref{22}, we get
\begin{equation}\label{ineq2}
\big\||x|^{-t-n/r}f\big\|_{L^r}\leq C\big\||x|^{-t+\alpha-n/s}(-\Delta)^{\alpha/2}f\big\|_{L^s}
\end{equation}
with the same $r,s$ in~\eqref{11},
since we are assuming $(-\Delta)^{\alpha/2}f\in C_0^\infty(\mathbb{R}^n\setminus\{0\})$.
Note that $1<q<n/\alpha$. So, we can choose $r_j$, $s_j$, $j=1,2$, such that
$$1<s_1<q<s_2<n/\alpha,\quad 1/s_j-1/r_j=\alpha/n,$$
and for $t_j=t+n(1/p-1/r_j)$
$$m-1<t_j<m,\quad\delta_t/2\leq\delta_{t_j}\leq3\delta_t/2.$$
Hence we can apply ~\eqref{ineq2} with $t=t_j$ to obtain
\begin{equation}\label{ineq3}
\big\||x|^{-t-n/p}f\big\|_{L^{r_j}}\leq C\big\||x|^{-t+\alpha-n/q}(-\Delta)^{\alpha/2}f\big\|_{L^{s_j}}
\end{equation}
for $j=1,2$.
Since $r_1<p<r_2$ and $s_1<q<s_2$, by real interpolation (\cite{St}) between the estimates in~\eqref{ineq3},
we see that for $1\leq w\leq\infty$
\begin{equation*}
\big\||x|^{-t-n/p}f\big\|_{L^{p,w}}\leq C\big\||x|^{-t+\alpha-n/q}(-\Delta)^{\alpha/2}f\big\|_{L^{q,w}}.
\end{equation*}
By choosing $w=q$, we get ~\eqref{ineq}.

Now we turn to the remaining case where $0<\alpha\leq n-1$.
In this case, ~\eqref{11} is valid for $1/s-1/r=\alpha/n$ and
\begin{equation}\label{55}
\frac12(1-\frac\alpha{n-1})+\frac\alpha n<\frac1s<\frac12+\frac\alpha{2(n-1)},
\end{equation}
because we can choose $\beta$ so that $n-1<\beta<n-\delta_t$.
Since ~\eqref{55} holds for $s$ replaced by $q$, repeating the previous argument,
one can show ~\eqref{ineq}. We omit the details.
\end{proof}



\begin{thebibliography}{9}

\bibitem{JK} D. Jerison and C. E. Kenig, \textit{Unique continuation and absence of positive eigenvalues for
Schr\"odinger operators}, Ann. of Math. 121 (1985), 463-494.

\bibitem{KN} C. E. Kenig and N. Nadirashvili, \textit{A counterexample in unique continuation},
Math. Res. Lett. 7 (2000), 625-630.

\bibitem{KT} H. Koch and D. Tataru, \textit{Sharp counterexamples in unique continuation for second order
elliptic equations}, J. Reine Angew. Math. 542 (2002), 133-146.

\bibitem{L} I. Laba, \textit{Unique continuation for Schr\"odinger operators and for higher powers of
the Laplacian}, Math. Methods Appl. Sci. 10 (1988), 531-542.

\bibitem{La} N. Laskin, \textit{Fractional quantum mechanics and L\'{e}vy path integrals},
Phys. Lett. A 268 (2000), 298-305.

\bibitem{S} I. Seo, \textit{Remark on unique continuation for higher powers of the Laplace operator},
J. Math. Anal. Appl. 397 (2013), 766-771.

\bibitem{S2} I. Seo, \textit{On unique continuation for Schr\"odinger operators of fractional and higher orders},
to appear in Math. Nachr., arXiv:1301.2460v2.

\bibitem{St} E. M. Stein, \textit{Interpolation of linear operators}, Tran. Amer. Math. Soc. 83 (1956), 482-492.   

\bibitem{St2} E. M. Stein, \textit{Appendix to ``unique continuation''}, Ann. of Math. 121 (1985), 489-494.

\end{thebibliography}
\end{document}